# HILBERT–PÓLYA CONJECTURE, ZETA–FUNCTIONS AND BOSONIC QUANTUM FIELD THEORIES

JULIO C. ANDRADE*
*Institute for Computational and Experimental Research in Mathematics (ICERM)*
*Brown University*
*121 South Main Street*
*Providence, RI 02903, USA.*
*julio_andrade@brown.edu*

*Dedicated to the 90th birthday of Professor Freeman Dyson.*

The original Hilbert and Pólya conjecture is the assertion that the non–trivial zeros of the Riemann zeta function can be the spectrum of a self–adjoint operator. So far no such operator was found. However the suggestion of Hilbert and Pólya, in the context of spectral theory, can be extended to approach other problems and so it is natural to ask if there is a quantum mechanical system related to other sequences of numbers which are originated and motivated by Number Theory.

In this paper we show that the functional integrals associated with a hypothetical class of physical systems described by self–adjoint operators associated with bosonic fields whose spectra is given by three different sequence of numbers cannot be constructed. The common feature of the sequence of numbers considered here, which causes the impossibility of zeta regularizations, is that the various Dirichlet series attached to such sequences - such as those which are sums over "primes" of (norm $P$)$^{-s}$ have a natural boundary, i.e., they cannot be continued beyond the line Re($s$) = 0. The main argument is that once the regularized determinant of a Laplacian is meromorphic in $s$, it follows that the series considered above cannot be a regularized determinant. In other words we show that the generating functional of connected Schwinger functions of the associated quantum field theories cannot be constructed.

## I. INTRODUCTION

The Riemann zeta function is the function of complex variable $s = \sigma + it$ defined to be the absolutely convergent series in the region Re($s$) > 1 by

$$\zeta(s) := \sum_{n=1}^{\infty} \frac{1}{n^s}, \qquad \mathrm{Re}(s) > 1. \qquad (1)$$

Riemann [1] proved that $\zeta(s)$ has an analytic continuation to the whole complex plane $\mathbb{C}$ except for simple pole at $s = 1$ with residue 1 at this pole.

In this same region, Re($s$) > 1, the Riemann zeta function can be represented as a product over all primes $p$, called Euler product, namely

$$\zeta(s) = \prod_{p \text{ prime}} \left(1 - \frac{1}{p^s}\right)^{-1}, \qquad (2)$$

---

* julio_andrade@brown.edu

and we can quickly see that the formula above connects properties of prime numbers with properties of the Riemann zeta function $\zeta(s)$.

One of the most acclaimed open problems in pure mathematics is the Riemann Hypothesis, which states that the nontrivial zeros of $\zeta(s)$, i.e., the zeros of $\zeta(s)$ on $0 \leq \sigma < 1$, have real part equal to $\frac{1}{2}$ (see Refs. [2–4]). The connection between the prime numbers and the zeros of Riemann zeta function can be seen from the assertion that the validity of the Riemann hypothesis is equivalent to the following asymptotic formula

$$\pi(x) = \mathrm{Li}(x) + O(\sqrt{x} \log x) \qquad \text{as } x \to \infty, \qquad (3)$$

where Li($x$) is the logarithmic integral function and $\pi(x)$ denotes the number of primes $p \leq x$.

The Hilbert–Pólya conjecture is the claim that the imaginary parts $t_n$ of the zeros $\frac{1}{2} + it_n$ of the Riemann zeta function are the eigenvalues of a self-adjoint (Hermitian) operator in a suitable Hilbert space and a proof of the Riemann hypothesis follows if we are able to find such an operator, but so far no such operator was found. Montgomery[5], inspired by an observation made by Freeman Dyson and later supported by numerical evidences from Odlyzko [6], conjectured that distribution of spacing between (pair–correlation) the Riemann zeros is the same as the pair correlation function of the eigenvalues of random matrices from the Gaussian unitary ensem-ble [7]. The results obtained by Montgomery on the pair–correlation of the zeros of the $\zeta(s)$ has lead some authors to consider and to investigate the Riemann hypothesis under the light of random matrix theory and quantum mechanics of classically chaotic systems [8, 9]. For a detailed discussion about the Riemann hypothesis and physics, see Refs. [8, 10–18]. On the assumption of the Hilbert–Pólya conjecture and the connection between zeros of $\zeta(s)$ and prime numbers it is natural to ask if there is a spectral interpretation for the sequence of prime numbers and other different number theoretical sequences of numbers. For example, there is a quantum mechanical system such that the Hamiltonian associated to it has the prime numbers as its spectrum?

Several authors [16, 19–24] have considered the question whether there is a quantum mechanical potential related to the prime numbers. In a recent and beautiful paper, Menezes and Svaiter [25] showed that for free scalar quantum field theories the associated functional integral cannot be constructed. The aim of this note is to extend the work of Menezes and Svaiter and show that there are different number theoretical sequences that not allow us to construct the generating functional of connected Schwinger functions of free scalar bosonic quantum field theories. The sequences we present in this note that causes the impossibility of construction of the functional integral of such theories are: the sequence of the primes $p \equiv 1 \pmod{m}$ with $m = 3, 4, 6$; the sequence given by $(N_{k/\mathbb{Q}}P)$, when $P$ varies over all prime ideal of $\mathcal{O}_k$ with $N_{k/\mathbb{Q}}P = p^{f_P}$ denoting the absolute norm of $P$ where $f_P = \deg P$ is the degree of the residue class field extension and $\mathcal{O}_k$ is the ring of integers of an algebraic number field $k$ such that $[k : \mathbb{Q}] = N$ (see Ref.[26]); and finally the sequence given by $(N\mathfrak{p})$ where $\mathfrak{p}$ runs over all prime divisors of the global function field $K/\mathbb{F}_q$ (see Section III for more details).

The common feature of the above sequences which causes the impossibility of zeta regularizations is that the various Dirichlet series attached to such sequences - such as those which are sums over "primes" of $(\text{norm } P)^{-s}$ have a natural boundary, i.e., they cannot be continued beyond the line $\text{Re}(s) = 0$. The main argument used in this paper is that once the regularized determinant of a Laplacian is meromorphic in $s$, it follows that the series considered above cannot be a regularized determinant. The situation is different if we consider, for example, the manifold as being the circle of dimension $N = 1$, then the eigenvalues of the Laplacian are given by $n^2$ for integers $n$. The Minakshisundaram-Pleijel zeta function is given by

$$Z(s) = \sum_{n \neq 0} \frac{1}{(n^2)^s} = \zeta(2s), \qquad (4)$$

and it is clear from properties of $\zeta(s)$ that $Z(s)$ has analytic continuation to the whole complex plane with a simple pole at $s = \frac{1}{2}$ and therefore does not have a natural boundary.

## II. QUANTUM FIELD THEORY AND SPECTRAL ZETA FUNCTION

We can say that, informally, Quantum Field Theory (QFT) is the extension of quantum mechanics (QM) in the sense that QFT handle with fields instead only particles and therefore the main aim of QFT is to deal with fields that are similar to a classical field, i.e., is a function defined over space and time, but which also incorporates the background of quantum mechanics. For example we can quote the quantum electrodynamics (QED), which is a quantum field theory and in essence is the relativistic quantum field theory of electrodynamics which has one electron field and one photon field and such theory was completed by Dyson [27] in 1949. The novelty involving QFT is that it handles with fields instead particle and so it is able to describe systems with an infinite number of degrees of freedom. The mathematical study of QFT is in part devoted to the properties of the resulting functional integral of the solutions of field theory models and make these functional integrals mathematically precise. By standard Euclidean methods it is possible to present a very well–defined functional integral representation, sometimes called partition function.

One of the main characteristic of the path integral formalism [28] for developing quantum field theories is that the spectral properties of the Hamiltonian of such theories can be deduced from properties of moment of the measure–theoretic factors in the functional integral. For example, the bosonic field theory in one–dimensional time can be reduced to the Heisenberg quantum mechanics formulation, i.e., given a quantum mechanical system with $n$ degrees of freedom we can find a $n$–component free scalar field model in one dimension which derives the desired quantum mechanical system.

The question if quantum field theory can be used to study number theory problems is not new and was considered by some authors [16, 20–24, 29, 30] and recently by Menezes and Svaiter [25]. In this note, as said earlier, we consider different number theoretical sequences of numbers and ask about the possibility of quantum field theories whose spectrum is given by such special number theoretical sequences.

We assume for simplicity that our quantum field is a neutral free scalar field $\varphi$ defined in a compact Riemannian $C^\infty$ manifold $(M, g)$ with metric $g = (g_{ij})$ on $M$. However, with extra care, the same results of this paper can be derived for higher spin fields. We also assume that $\hbar = c = 1$.

From now on we consider a free neutral scalar field $\varphi$ defined on $M$ and we assume $M$ is a $d$–dimensional Minkowski spacetime ($\dim M = d$) and $M$ is connected. The Euclidean field theory is obtained by analytic continuation of the $n$–point Wightman functions to imaginary time. Let us assume that the Euclidean space is compact and may or may not have boundary. If the Euclidean space has a smooth boundary we can find appropriate boundary conditions on the fields. Also let us assume there exists an elliptic, self–adjoit differential operator $D$ which acts on scalar functions $\varphi$ on the Euclidean space.

Consider one of the most common examples for $D$, which is the negative of the Laplacian shifted by the mass of the model, given by

$$D = (-\Delta + m_0^2), \qquad (5)$$

where $\Delta$ is the $d$–dimensional Laplacian and $m_0$ is the mass of the model. Denote the kernel by

$$K(m_0; x - y) = (-\Delta + m_0^2)\delta^d(x - y), \qquad (6)$$

hence in this case the Euclidean generating functional $Z[h]$ is defined by the following functional integral

$$Z[h] = \int [d\varphi] \exp\left(-S_0 + \int d^d x h(x)\sqrt{g(x)}\varphi(x)\right), \qquad (7)$$

where $g = \det(g_{ij})$ and $[d\varphi]$ is a translational invariant measure given by

$$[d\varphi] = \prod_x d\varphi(x), \qquad (8)$$

with $h(x)$ being a smooth function generating the Schwinger functions of the theory and $S_0$ is the action which describes a free scalar field given by

$$S_0(\varphi) = \int d^d x d^d y \sqrt{g(x)g(y)}\varphi(x)K(m_0; x - y)\varphi(y). \qquad (9)$$

To generate the connected Schwinger functions of the theory we define the functional $W[h] = \ln Z[h]$ which will generate such functions. Hawking [31] showed that in flat space, in which the eigenvalues of Laplacians are known, the zeta function corresponding to the partition function can be computed explicitly and he studied zeta function regularization in order to calculate the partition functions for thermal graviton and matter's quanta in the horizon of black holes and on de Sitter space using the relation by the inverse Mellin transformation to the trace of the kernel of heat equations. In our case we also need to performing the zeta function regularization, i.e., we need to regularize a determinant associated with the operator $D$, once $W[0] = -\frac{1}{2}\ln \det D$, and thus we will be able to obtain a well–defined object.

The spectral zeta function is defined to be

$$\zeta_A(s) := \sum_{n \in \mathbb{N}} \frac{1}{a_n^s} \qquad \text{Re}(s) > s_0, \qquad (10)$$

for some complex number $s_0$ and where $A = (a_n)_{n \in \mathbb{N}}$ is a sequence of nonzero complex numbers. Provided that $\zeta_A(s)$ has analytic extension and is holomorphic at $s = 0$ the regularization zeta function product is defined as

$$\prod_{n \in \mathbb{N}} a_n := \exp\left(-\frac{d}{ds}\zeta_A(s)\bigg|_{s=0}\right). \qquad (11)$$

Using the standard machinery of spectral theory of elliptic operators [32] we know that there exists a complete orthonormal set $\{f_k\}_{k=1}^{\infty}$ such that the associated eigenvalues obey: $0 \leq \lambda_1 \leq \lambda_2 \leq \cdots \leq \lambda_k \to \infty$ when $k \to \infty$ (The zero eigenvalue $\lambda_0$ must be omitted and the eigenvalues are counted with multiplicities). From Linear Algebra we get that the operator $D$ in the basis $\{f_k\}$ is represented by the infinite diagonal matrix

$$D = \begin{pmatrix} \lambda_1 & & & \\ & \lambda_2 & & \\ & & \lambda_3 & \\ & & & \ddots \end{pmatrix}, \qquad (12)$$

and therefore the $f_n$'s satisfy the eigenvalue equation for $D$, namely

$$Df_n(x) = \lambda_n f_n(x). \qquad (13)$$

The spectral zeta function $\zeta_D(s)$ attached to the operator $D$ is defined to be the infinite series

$$\zeta_D(s) = \sum_{n \in \mathbb{N}} \frac{1}{\lambda_n^s}, \qquad \text{Re}(s) > s_0, \qquad (14)$$

for $s_0$ sufficiently large and using the ideas above we have, formally, that

$$-\frac{d}{ds}\zeta_D(s)\bigg|_{s=0} = \ln \det D. \qquad (15)$$

The spectral zeta function $\zeta_D(s)$ must be defined at $s = 0$ as can be seen from (15). In particular to be able to regularize the determinant of the functional integral of the theory it is necessary that $\zeta_D(s)$ has an analytic continuation into the whole complex plane. We notice that the spectral zeta function satisfies the following scaling property

$$\frac{d}{ds}\zeta_{\mu^2 D}(s)\bigg|_{s=0} = \ln \mu^2 \zeta_D(s)\bigg|_{s=0} + \frac{d}{ds}\zeta_D(s)\bigg|_{s=0}, \qquad (16)$$

where $\mu$ is an arbitrary parameter with dimension of mass.

## III. ZETA FUNCTION OF CURVES AND PRIME DIVISORS

We start by defining the zeta function of a curve $C$. Let $C$ be a nonsingular projective curve over a finite field $\mathbb{F}_q$

of characteristic $p$ with $q = p^a$ where $p$ is a prime number (The prototype for such finite fields is $\mathbb{Z}/p\mathbb{Z}$). Let $\mathcal{D}_C$ be the usual additive group of divisor on $C$ defined over $\mathbb{F}_q$. In other words they are finite sums $\mathfrak{a} = \sum a_i P_i$ where $a_i \in \mathbb{Z}$ and the $P_i$'s are points of $C$ defined over a finite extension of $\mathbb{F}_q$, such that $\phi(\mathfrak{a}) = \mathfrak{a}$ where $\phi$ is the Frobenius endomorphism on $C$ raising the coordinates to the $q$th power. We call $\deg(\mathfrak{a}) = \sum a_i$ the degree of the divisor $\mathfrak{a}$ and it is called effective if every $a_i > 0$ and in this case we write $\mathfrak{a} > 0$.

A prime divisor $\mathfrak{p}$ is a positive divisor that cannot be expressed as the sum of positive divisors, in other words they are the analogue of the prime numbers for the function field of the curve $C$. The norm of a divisor $\mathfrak{a}$ is defined as $N\mathfrak{a} = q^{\deg(\mathfrak{a})}$. With this notation we define the zeta function of the curve $C$ to be

$$\zeta(s, C) = \sum_{\mathfrak{a} > 0} \frac{1}{N\mathfrak{a}^s} \qquad \mathrm{Re}(s) > 1, \qquad (17)$$

(for more details see [33, Chapter 5]).

Making use of the multiplicativity of the norm and the fact that the group of divisors of $C$, $\mathcal{D}_C$, is a free abelian group on the set of prime divisors, we see that the zeta–function of $C$ also has an Euler product

$$\zeta(s, C) = \prod_{\mathfrak{p}} (1 - N\mathfrak{p}^{-s})^{-1} \qquad \mathrm{Re}(s) > 1, \qquad (18)$$

and the Poincaré duality gives us the functional equation

$$q^{(g-1)s}\zeta(s, C) = q^{(g-1)(1-s)}\zeta(1-s, C), \qquad (19)$$

where $g$ is the genus of the curve $C$. It is known [33] that (17) has an analytic continuation to all $\mathbb{C}$ and has simple poles at $s = 0$ and $s = 1$. Weil [34] showed that the Riemann hypothesis is true for $\zeta(s, C)$, i.e., all the zeros of $\zeta(s, C)$ lie on the line $\mathrm{Re}(s) = 1/2$.

Let us assume that there is an hypothetical operator $D$, given as in the Section II, which acts on scalar functions in the Riemannian manifold $M$ and has the sequence of $(N\mathfrak{p})$ as its spectrum when $\mathfrak{p}$ varies over the prime divisors of $C$. Therefore the spectral zeta function, $\zeta_D(s)$, attached to $D$ is given by

$$P(s, C) = \sum_{\mathfrak{p}} \frac{1}{N\mathfrak{p}^s} \qquad \mathrm{Re}(s) > 1, \qquad (20)$$

where the sum is over all prime divisors on $C$. Using the same reasoning from Section II we know, in this case, that to be able to regularize the generating functional of connected Schwinger functions of the free scalar quantum field theory we must analytically extend $P(s, C)$ to $s = 0$. So now our main aim is to study the analytic continuation of (20).

Using (18) we have

$$\log \zeta(s, C) = \sum_{n=1}^{\infty} \frac{1}{n} P(ns, C), \qquad \mathrm{Re}(s) > 1. \qquad (21)$$

Invoking the Möbius inversion formula, one has

$$P(s, C) = \sum_{k=1}^{\infty} \frac{\mu(k)}{k} \log \zeta(ks, C) \qquad \mathrm{Re}(s) > 1, \qquad (22)$$

where the Möbius function $\mu(n)$ is defined for $n > 1$ as [35]

$$\mu(n) = \begin{cases} (-1)^k & \text{if } n = p_1 \cdots p_k \text{ for distinct primes} \\ 0 & \text{otherwise} \end{cases},$$

and $\mu(1) = 1$. From (22) we have that the analytic continuation of (20) can be obtained from the analytic continuation of the zeta function $\zeta(s, C)$ of the curve $C$. The analytic continuation of $\zeta(s, C)$ was proved by Schmidt and Weil and we know from them that $\zeta(s, C)$ has an analytic continuation to all $\mathbb{C}$ with simple poles at $s = 0$ and $s = 1$ and the analytic continuation is given by

$$\zeta(s, C) = \frac{L_C(q^{-s})}{(1 - q^{-s})(1 - q^{1-s})}, \qquad (23)$$

where $L_C(q^{-s})$ is a polynomial in $u = q^{-s}$ of degree $2g$ with integer coefficients.

For $\mathrm{Re}(s) \leq 1$ the analytic structure of $P(s, C)$ is given by the poles and zeros of $\zeta(s, C)$ and from (22) one notices that $s = \frac{1}{k}$ is a singularity of $P(s, C)$ for all square–free positive integers $k$. We have that the sequence $\left(\frac{1}{k}\right)$ limits to $s = 0$.

We have that all points on the line $\mathrm{Re}(s) = 0$ are limit points of the poles of $P(s, C)$ and so we can conclude that the line $\mathrm{Re}(s) = 0$ is a natural boundary of $P(s, C)$. Indeed, the same arguments presented in Refs. [36, 37] for the Riemann zeta function $\zeta(s)$ holds for the $\zeta(s, C)$ and hence there is a clustering of singular points along the imaginary axis from the zeros of $\zeta(s, C)$ and we can see $\mathrm{Re}(s) = 0$ is a natural boundary for $P(s, C)$. Therefore we can conclude that the function $P(s, C)$ can be analytically continued only in the strip $0 < \sigma \leq 1$, and with the exception of the singular points of $P(s, C)$ in $0 < \sigma \leq 1$ the equations (22) and (23), gives a representation for $P(s, C)$ which is valid for $\mathrm{Re}(s) > 0$.

Denoting $\zeta(ks) = z$, where $z$ is a complex variable and using that the function $\ln z$ is analytic at every point of its Riemann surface and satisfies

$$\frac{d}{dz} \ln z = \frac{1}{z}, \qquad (24)$$

we obtain that $dP(s, C)/ds = P'(s)$ can be computed for $\mathrm{Re}(s) > 0$. From (22), one has

$$\frac{d}{ds}P(s,C) = \sum_{k=1}^{\infty} \frac{\mu(k)}{k} \frac{1}{\zeta(ks,C)} \frac{\partial}{\partial s}\zeta(ks,C), \quad \mathrm{Re}(s) > 0. \tag{25}$$

In order to obtain a well–defined functional integral associated to the scalar field $\varphi$, as mentioned in Section II, we need to define the regularized determinant associated to $D$, i.e., it is necessary to compute the corresponding spectral zeta function and its derivative at $s=0$. For the spectrum given by the sequence of numbers $(N\mathfrak{p})$, the functional $W[0]$ is given by

$$W[0] = \frac{1}{2}\ln\mu^2 P(s,C)\big|_{s=0} + \frac{1}{2}\frac{d}{ds}P(s,C)\big|_{s=0}, \tag{26}$$

where we have used that $W[0] = -\frac{1}{2}\ln\det D$, and equations (15) and (16).

Combining (25) and (26) we have that if $k$ is a square–free number in $\zeta(ks,C) = 0$ then $dP(s,C)/ds$ diverges. And if $k$ is not a square–free number the quantity $dP(s,C)/ds$ is not well–defined. Hence the sequence given by $(N\mathfrak{p})$ of the general curve $C$ leads us to a ill–defined functional integral for the quantum field theory in question. In other words we have proved the following theorem:

**Theorem 1** *The sequence of numbers given by $(N\mathfrak{p})$ where $\mathfrak{p}$ runs over all prime divisors of the curve $C$ is not zeta regularizable, i.e., there is no hypothetical free scalar bosonic quantum field theories described by self–adjoint operators $D$ such that its spectrum is given by $(N\mathfrak{p})$.*

Now we will show that the sequence of the primes $p \equiv 1 (\mathrm{mod}\ m)$ with $m = 3, 4, 6$ and the sequence given by $(N_{k/\mathbb{Q}}P)$, when $P$ varies over all prime ideal of $\mathcal{O}_k$ with $N_{k/\mathbb{Q}}P = p^{f_P}$ denoting the absolute norm of $P$ where $f_P = \deg P$ is the degree of the residue class field extension and $\mathcal{O}_k$ is the ring of integers of an algebraic number field $k$ such that $[k:\mathbb{Q}] = N$, are not zeta regularizable.

For the sequence of primes $p \equiv 1(\mathrm{mod}\ m)$ with $m = 3, 4, 6$, Kurokawa [38] showed that

$$\sum_{p \equiv 1(\mathrm{mod}\ m)} \frac{1}{p^s}, \tag{27}$$

has analytic continuation for $\mathrm{Re}(s) > 0$ with the natural boundary $\mathrm{Re}(s) = 0$. Hence the same arguments above apply to this case. And for the sequence given by $(N_{k/\mathbb{Q}}P)$, analogous calculations such as those carried out in Section III can be performed for $(N_{k/\mathbb{Q}}P)$, where in this case, we use the Dedekind zeta function of the algebraic function field $k$. Therefore this imply that the sequences above lead us to ill–defined free scalar quantum field theories.

## IV. CONCLUSION

We have showed in this paper that, beyond the sequence of prime numbers, other number theoretical sequences of numbers lead us to ill–defined free neutral scalar bosonic quantum field theories described by self-adjoint operators $D$ acting on scalar fields $\varphi$ in a $d$–dimensional Minkowski spacetime.

In two recent works, by Sierra and Townsend [17] and Sierra and Rodriguez-Laguna [18], the Berry–Keating classical Hamiltonian $H = xp$, conjectured to be the semi–classical limit of a quantum mechanical model for the complex zeros of the Riemann zeta function, is revisited. In the first work they showed that classical results can be recovered by viewing the $H = xp$ model as a lowest Landau level limit of a quantum mechanical model for a charged particle on the $xy$–plane and so the authors provide semi-classical evidence that the full 'counting' formula for the Riemann zeros will arise from a consideration of the higher Landau levels. In the second work, the authors reformulate the Berry-Keating model in terms of the classical Hamiltonian $H_{\mathrm{cl}} = x(p+l_p^2/p)$ defined on the half-line $x \geq l_x$ and in this case it has closed orbits whose semiclassical spectrum agrees with the average Riemann zeros. To achieve the goal of finding the Riemann operator we need to understand the quantum origin of the fluctuations of the Riemann zeros and the work of Sierra and Rodriguez–Laguna [18] suggest that is needed to modify the Hamiltonian in a nontrivial way. Compared with our case, the operator $D = (-\Delta + m_0^2)$ is the Laplacian on a $d$–dimensional Riemannian manifold acting on a free scalar field instead only a particle. And for the cases considered above, the spectrum of $D$, lead us to ill–defined free scalar quantum field theories. It would be interesting to find a classical Hamiltonian, such as those suggested by Berry-Keating or the Sierra–Rodriguez–Laguna that would provide us the prime numbers and other number theoretical sequences as its spectrum and at the same time lead us to well–defined quantum physical systems.

As we have remarked before, the same results above can be extended to higher spin fields, in particular spinor fields, with weak interactions and also for quantum mechanical systems with $n$ degrees of freedom since they can be derived from $n$-component scalar quantum field theories in one dimension. The results obtained in this paper are viewed as constraints for formulations of quantum mechanical systems which has the sequences of numbers as presented here as its spectrum. And if there is a hypothetical quantum system which has an energy spectrum given by those sequences of the theorem above this will imply that cannot be a free energy associated with such a system. But an intriguing issue is whether this situation persists in fermionic quantum field theories and such question is under investigation by the author.

Other very interesting question to be investigated is if there are other sequences of prime numbers, such as the twin primes which obey the same result above. Now, returning to the case of zeta functions of curves $C$, we know that there is a spectral interpretation for the zeros of the $\zeta(s,C)$ through the Frobenius on Cohomology. The question if there is a quantum field theory or a physical interpretation for the Frobenius can shed some light in the path of a proof of the original Riemann Hypothesis and in this way a physical interpretation for the Frobenius using the ideas of this paper will be very interesting and is being carried out by the same author.

## V. ACKNOWLEDGMENTS

The author would like to thank professors Freeman Dyson, Michael Rosen and Peter Sarnak for the encouragement and several useful discussions. We also would like to thank the valuable comments, suggestions and a theoretical clarification from an anonymous referee which helped to greatly improve the presentation of this paper.

The author was supported by a NFS postdoctoral research fellowship and an ICERM–Brown University Postdoctoral Research Fellowship.